\documentclass[11pt,twoside]{article}
\usepackage[left=3cm,right=3cm,top=3cm,bottom=3cm]{geometry}
\usepackage{amssymb}
\usepackage{graphicx}
\usepackage{graphicx,type1cm,eso-pic,color}
\usepackage{amssymb}
\usepackage{amssymb,amsmath}
\usepackage{graphicx,epsfig}
\usepackage{enumerate}
\usepackage{color}
\usepackage {multicol}

\numberwithin{equation}{section}
\newtheorem{theorem}{Theorem}
\newtheorem{lemma}{Lemma}

\newtheorem{proposition}{Proposition}
\newtheorem{remark}{Remark}
\newtheorem{corollary}{Corollary}
\newtheorem{example}{Example}
\makeatletter
\AddToShipoutPicture{
	\setlength{\@tempdimb}{.5\paperwidth}
	\setlength{\@tempdimc}{.5\paperheight}
	\setlength{\unitlength}{1pt}
	\put(\strip@pt\@tempdimb,\strip@pt\@tempdimc)
}
\makeatother

\usepackage{scalerel,stackengine}
\stackMath
\newcommand\reallywidehat[1]{%
	\savestack{\tmpbox}{\stretchto{%
			\scaleto{%
				\scalerel*[\widthof{\ensuremath{#1}}]{\kern-.6pt\bigwedge\kern-.6pt}%
				{\rule[-\textheight/2]{1ex}{\textheight}}
			}{\textheight}%
		}{0.5ex}}%
	\stackon[1pt]{#1}{\tmpbox}%
}
\begin{document}
	\setcounter{page}{1}

	\thispagestyle{empty}
	\markboth{}{}

	\pagestyle{myheadings}
	\markboth{ S.K.Chaudhary, N. Gupta }{ S.K.Chaudhary, N. Gupta }
	
	\date{}
	
	
	\noindent  
	
	\vspace{.1in}
	
	{\baselineskip 20truept
		
		\begin{center}
			{\Large {\bf General weighted extropy of minimum and maximum ranked set sampling with unequal samples}} \footnote{\noindent	{\bf * } Corresponding author E-mail: skchaudhary1994@kgpian.iitkgp.ac.in\\
				{\bf **} E-mail: nitin.gupta@maths.iitkgp.ac.in
			}			
		\end{center}

		\vspace{.1in}
		
		\begin{center}
			{\large {\bf Santosh Kumar Chaudhary*, Nitin Gupta**}}\\
			{\large {\it Department of Mathematics, Indian Institute of Technology Kharagpur, West Bengal 721302, India }}
			\\
		\end{center}

		\vspace{.1in}
		\baselineskip 12truept

		\begin{abstract}
			In industrial, environmental, and ecological investigations, ranked set sampling is a sample method that enables the experimenter to use the whole range of population values. The ranked set sampling process can be modified in two extremely helpful ways: maximum ranked set sampling with unequal samples and minimum ranked set sampling with unequal samples. They permit an increase in set size without too many ranking errors being introduced. In this paper, we are defining general weighted extropy (GWJ) of minimum and maximum ranked set samples when samples are of unequal size (minRSSU and maxRSSU, respectively). Stochastic comparison and monotone properties have been studied under different situations. Additionally, we compare the extropy of these two sampling data with that of ranked set sampling data and simple random sampling data. Finally, Bounds of GWJ of minRSSU and maxRSSU have been obtained.

			\vspace{.1in}
			
			
		\end{abstract}
		
		\section{Introduction}\label{section1}
		McIntyre (1952) suggested ranked set sampling (RSS) as a way to improve the sample mean's accuracy as an estimate of the population means. In terms of estimating the population mean, it appears that the RSS is a more effective sampling technique than the simple random sample (SRS). A sampling technique called RSS can be used to increase the cost-effectiveness of choosing sample units for an experiment or study. It is generally advised that ranking sample units instead of measuring them would be more affordable and straightforward. An example of the application of RSS as proposed by Stokes (1977) is given in Bain and Engelhardt (2017), where the study variate Y represents the oil pollution of sea water and auxiliary variable X represents the tar deposit in the nearby seashore. Here, collecting seawater samples and measuring oil pollution is difficult and costly. However, the prevalence of pollution in seawater is much reflected in the tar deposit in the surrounding terminal seashore. In this example, ranking the pollution level of seawater based on the tar deposit on the seashore is more natural and scientific than ranking it visually or by judgment method. In comparison to SRS, the RSS is a more effective sampling method for estimating the population mean. Readers can refer to Ozturk and Wolfe (2000), Chen et al. (2004), Jemain et al.(2007), Tseng and Wu (2007), and Al-Saleh and Diab (2009) for more details on RSS.		
		
		Al-Odat and Al-Saleh (2001) and Biradar and Santosha (2014) presented minimal ranked set sampling with unequal samples (MinRSSU) and maximal ranked set sampling with unequal samples (maxRSSU), respectively, as two significant RSS variants. They viewed the MinRSSU and maxRSSU as moving extreme ranked set sampling designs, however, Rahmani and Razmkhah (2017) looked into those individually to achieve the ideal ranking test. As a helpful improvement to the RSS technique, Qiu and Eftekharian (2021) investigated the extropy of the minimum and maximum ranked set sampling procedure with unequal samples (minRSSU and maxRSSU). In the minRSSU and maxRSSU, we draw $n$ simple random samples, in which the size of the $i-$th samples is $i, \  i = 1,\dots n.$ The one-cycle minRSSU or maxRSSU involves an initial ranking of $n$ samples of size $n$.  We collect $Z_i=X_{(1:i)i}$ for all $i=1,2,\dots, n, \ $ and $Y_i=X_{(i:i)i}$ for all $i=1,2,\dots, n, \ $ where $X_{(i:j)j}$ denotes the $i$-th order statistic from the $j$-th SRS of size $j.$ The resulting sample is called one-cycle minRSSU and maxRSSU of size $n$ and denoted by $Z_{minRSSU} =\{Z_i : i=1,\dots,n\}$ and $Y_{maxRSSU} =\{Y_i : i=1,\dots,n\},$ respectively (see Biradar and Santosha (2014)).
		
		Measures of information for RSS and its variations have been developed by several authors (see Eskandarzadeh et al.(2018), Raqab and Qiu (2019), Qiu and Raqab (2022) and Gupta and Chaudhary (2023) etc). Let X be a non-negative absolutely continuous random variable with probability density function (pdf) $f.$ Lad et al. (2015) developed a new information measure called extropy,
		\begin{equation}\label{extropy}
			J(X)=-\frac{1}{2} \int_{0}^{\infty}f^2(x)dx=-\frac{1}{2}E\left(f(X)\right).
		\end{equation}		
		Extropy was used by Raqab and Qiu (2019) for evaluating the information content of RSS data. Qiu and Eftekharian (2021) considered the information content of minimum and maximum ranked set sampling with unequal samples in terms of extropy. The extropy properties of the ranked set sample (RSS) when the ranking is not perfect are considered by Chacko and George (2022). By deriving the expression for extropy of concomitant order statistic, the expression for extropy of RSS of the study variable Y in which an auxiliary variable $X$ is used to rank the units in each set, under the assumption that $(X, Y)$ follows Morgenstern family of distributions is obtained by Chacko and George (2022). 
		
		Gupta and Chaudhary (2023) introduced the  general weighted extropy (GWJ) for the non-negative random variable $X$ with weight $w(x)\geq 0$ as
		\begin{align}
			J^{w}(X)&=-\frac{1}{2} \int_{0}^{\infty}w(x) f^2(x)dx =-\frac{1}{2} \int_{0}^{\infty} w(F^{-1}(u))f(F^{-1}(u))du=-\frac{1}{2}E(\Lambda_X^{w}(U)),
		\end{align}
		where $\Lambda_X^{w}(u)=w(F^{-1}(u))f(F^{-1}(u))$ and $U$ is uniformly distributed random variable on $(0,1)$, i.e., $U\sim$ uniform$(0,1)$.

		This paper's main objective is to study the general weighted extropy characteristics of maxRSSU and minRSSU data. The rest of this paper is organized as follows. Section 2 defines the GWJ of minRSSU and maxRSSU and evaluated for different distributions.  In Section 3, stochastic comparisons have been done. Section 4 is devoted to the monotone properties of GWJ of minRSSU and maxRSSU. Bounds for GWJ of minRSSU and maxRSSU have been given in section 5. Section 6 concludes this paper.

		\section{GWJ of maxRSSU and minRSSU}
		Let $X$ be a non-negative random variable with finite mean $\mu$ and variance $\sigma^2$. Let $X_1,X_2,\dots,X_n$ be a random sample from $X.$ The joint pdf for $\textbf{X}_{SRS}=\{X_i:\ i=1,\ldots,n\}$ is $\prod_{i=1}^{n}f(x_i).$ The general weighted extropy of $\textbf{X}_{SRS}^{(n)}$ is defined as (Gupta and Chaudhary (2023))
		\begin{align}
			J^{w}(\textbf{X}_{SRS}^{(n)})=\frac{-1}{2}\prod_{i=1}^{n}\left(\int_{0}^{\infty}w(x_i)f^2(x_i)dx_i\right) =\frac{-1}{2}\left(-2J^{w}(X)\right)^n =\frac{-1}{2}\left(E(\Lambda_X^{w}(U))\right)^n.
		\end{align}
		
		The GWJ of   $\textbf{X}_{RSS}^{(n)}$  has been given as (see Gupta and Chaudhary (2023))
		\begin{align} J^{w}(\textbf{X}_{RSS}^{(n)})=-\frac{1}{2}\prod_{i=1}^{n}\left(-2J^{w}(X_{(i:n)i})\right)\nonumber =- \frac{Q_n}{2}\prod_{i=1}^{n}E\left(\Lambda_X^{w} (B_{2i-1:2n-1})\right),
		\end{align}
		where \begin{align*}
			Q_n&=n^n\prod_{i=1}^{n}c_{i,n},\\
			c_{i,n}&=\frac{\binom{2i-2}{i-1}\binom{2n-2i}{n-i}}{\binom{2n-1}{n-1}}, 
		\end{align*}
		and $B_{2i-1:2n-1}$ is a beta distributed random variable with parameters $(2i-1)$ and $(2n-2i+1)$. 
		
		\noindent	Now, we can write the GWJ of   $\textbf{X}_{minRSSU}^{(n)}$  as
		\begin{align} J^{w}(\textbf{X}_{minRSSU}^{(n)})&=-\frac{1}{2}\prod_{i=1}^{n}[-2J^w(X_{1:i})] \label{minrssuformulaline1}\\ &=-\frac{1}{2}\prod_{i=1}^{n}\left(\int_{0}^{\infty}w(x)f^2_{1:i}(x)dx\right)\nonumber \\
			&=-\frac{1}{2}\prod_{i=1}^{n}\left(\int_{0}^{\infty}w(x) i^2 (1-F(x))^{2i-2}f^2(x)dx\right)\nonumber \\
			&=-\frac{(n!)^2}{2}\prod_{i=1}^{n}\left(\int_{0}^{1}w(F^{-1}(u)) (1-u)^{2i-2}f(F^{-1}(u))du\right)\nonumber \\
			&=-\frac{(n!)^2}{2(2n-1)!!}\prod_{i=1}^{n}\left(\int_{0}^{1} (2i-1)(1-u)^{2i-2}w(F^{-1}(u)) f(F^{-1}(u))du\right)\nonumber \\
			&=-\frac{(n!)^2}{2(2n-1)!!}\prod_{i=1}^{n} E(\Lambda_X^{w}(B_{1,2i-1}))\label{minrssuformula} ,
		\end{align}
		where $(2n-1)!!=\prod_{i=1}^{n}(2i-1)=1.3.5.7 \dots (2n-1)$ and $B_{1,2i-1}$ is a random variable having a beta distribution with parameters 1 and $2i-1.$

		\noindent	Also, we can write the GWJ of   $\textbf{X}_{maxRSSU}^{(n)}$  as
		\begin{align} J^{w}(\textbf{X}_{maxRSSU}^{(n)})&=-\frac{1}{2}\prod_{i=1}^{n}[-2J^w(X_{i:i})] \label{maxrssuformulaline1} \\&=-\frac{1}{2}\prod_{i=1}^{n}\left(\int_{0}^{\infty}w(x)f^2_{i:i}(x)dx\right)\nonumber \\
			&=-\frac{1}{2}\prod_{i=1}^{n}\left(\int_{0}^{\infty}w(x) i^2 (F(x))^{2i-2}f^2(x)dx\right)\nonumber \\
			&=-\frac{(n!)^2}{2}\prod_{i=1}^{n}\left(\int_{0}^{1}w(F^{-1}(u)) u^{2i-2}f(F^{-1}(u))du\right)\nonumber \\
			&=-\frac{(n!)^2}{2(2n-1)!!}\prod_{i=1}^{n}\left(\int_{0}^{1} (2i-1)u^{2i-2}w(F^{-1}(u)) f(F^{-1}(u))du\right)\nonumber \\
			&=-\frac{(n!)^2}{2(2n-1)!!}\prod_{i=1}^{n} E(\Lambda_X^{w}(B_{2i-1,1})) \label{maxrssuformula},
		\end{align}
		where $B_{2i-1,1}$ is a random variable having beta distribution with parameters $2i-1$ and 1.
		
		\begin{remark}
			For $w(x)=1,$ we obtain expressions given by Qiu and Eftekharian (2021).
		\end{remark}
		Now we provide some examples to illustrate the equation $(\ref{minrssuformula})$ and $(\ref{maxrssuformula})$.
		
		\begin{example}
			Let $Y$ be a random variable with power distribution. The pdf and cdf of $Y$ are respectively $f(y)= \theta y^{\theta -1}$ 
			and $F(y)= y^{\theta}$ , $0<y<1$ ,  $\theta > 0$. Let $w(y)=y^m,\ y>0, \ m>0$, then
			it follows that 
			\[\Lambda_Y^{w}(u)= w(F_Y^{-1}(u)) f_Y(F_Y^{-1}(u))=\theta u^{\frac{m + \theta - 1}{\theta}}\]
			for $ w(x)=x^m.$ 
			Then from (\ref{minrssuformula}) we have 
			\begin{align*}
				J^{w}(\textbf{X}_{minRSSU}^{(n)}) &=-\frac{(n!)^2}{2}\prod_{i=1}^{n}\left(\int_{0}^{1} (1-u)^{2i-2} w(F^{-1}(u)) f(F^{-1}(u))du\right)\nonumber \\
				&=-\frac{(n!)^2}{2}\prod_{i=1}^{n}\left(\int_{0}^{1} (1-u)^{2i-2} \theta u^{\frac{m + \theta - 1}{\theta}}\ du\right)\nonumber \\
				&=-\frac{(n!)^2 \theta^n}{2}\prod_{i=1}^{n}\left(\int_{0}^{1}  u^{\frac{m + \theta - 1}{\theta}}\ (1-u)^{2i-2} du\right)\nonumber \\
				&=-\frac{(n!)^2 \theta^n}{2}\prod_{i=1}^{n} \beta\left(\frac{m + \theta - 1}{\theta}+1, 2i-1\right),
			\end{align*}
			where $\beta(a,b)=\int_{0}^{1} u^{a-1}(1-u)^{b-1}du.$
			
			Also, from (\ref{maxrssuformula}) we have
			\begin{align*}
				J^{w}(\textbf{X}_{maxRSSU}^{(n)}) &=-\frac{(n!)^2}{2}\prod_{i=1}^{n}\left(\int_{0}^{1} u^{2i-2} w(F^{-1}(u)) f(F^{-1}(u))du\right)\nonumber \\
				&=-\frac{(n!)^2}{2}\prod_{i=1}^{n}\left(\int_{0}^{1} u^{2i-2} \theta u^{\frac{m + \theta - 1}{\theta}}du\right)\nonumber \\	
				&=-\frac{(n!)^2 \theta^n}{2}\prod_{i=1}^{n}\left(\int_{0}^{1} u^{2i-2+\frac{m + \theta - 1}{\theta}}du\right)\nonumber \\	
				&=-\frac{(n!)^2 \theta^n}{2}\prod_{i=1}^{n} \frac{1}{2i-1+\frac{m + \theta - 1}{\theta}}\nonumber \\	
				&=-\frac{(n!)^2 \theta^n}{2}\prod_{i=1}^{n}\frac{\theta}{2i\theta+m-1}.	
			\end{align*}
		\end{example}

		\begin{example}\label{eg2exp}
			Let $Z$ have an exponential distribution with cdf $F_Z(z)=1-e^{-\lambda z}, \ \lambda >0, \ z>0$. Let $w(x)=x^m, \ m>0, \ x>0$, then it follows that
			\begin{align*}
				\Lambda_Z^{w} (u) = w(F^{-1}(u)) f(F^{-1}(u))
				= \frac{(-1)^m (1-u) (ln(1-u))^m}{\lambda^{m-1}}, 0<u<1.  	
			\end{align*}
			Then from (\ref{minrssuformula}) we have 
			\begin{align*}
				J^{w}(\textbf{X}_{minRSSU}^{(n)}) &=-\frac{(n!)^2}{2}\prod_{i=1}^{n}\left(\int_{0}^{1} (1-u)^{2i-2} w(F^{-1}(u)) f(F^{-1}(u))du\right)\nonumber \\
				&=-\frac{(n!)^2}{2}\prod_{i=1}^{n}\left(\int_{0}^{1}  \frac{(1-u)^{2i-1}(-1)^m (ln(1-u))^m}{\lambda^{m-1}}du\right)\nonumber \\
				&=-\frac{(n!)^2}{2}\prod_{i=1}^{n}\frac{(-1)^m}{\lambda^{m-1}}\left(\int_{0}^{1} {(1-u)^{2i-1} (ln(1-u))^m}du\right).\nonumber \\
			\end{align*}
			Taking $u=1-e^{-x}$ and using \ $\int_{0}^{\infty} x^{\alpha-1} e^{-\beta x}dx=\frac{\Gamma(\alpha)}{\beta^\alpha}$ we obtain,
			\begin{align*}
				J^{w}(\textbf{X}_{minRSSU}^{(n)}) 
				&=-\frac{(n!)^2}{2}\prod_{i=1}^{n}\frac{(-1)^m}{\lambda^{m-1}}\left(\int_{0}^{1} {e^{-2ix} (-x)^m}dx\right).\nonumber \\
				&=-\frac{(n!)^2}{2}\prod_{i=1}^{n}\frac{1}{\lambda^{m-1}}\left(\int_{0}^{1} {e^{-2ix} x^m}dx\right).\nonumber \\
				&=-\frac{(n!)^2}{2}\prod_{i=1}^{n}\frac{1}{\lambda^{m-1}}\left(\frac{\Gamma(m+1)}{(2i)^{m+1}}\right).\nonumber \\
				&=-\frac{1}{2}\left(\frac{1}{\lambda^{m-1}}\right)^{n}\left(\frac{\Gamma(m+1)}{2^{m+1}}\right)^n \left(\frac{1}{n!}\right)^{m-1}.\nonumber \\
			\end{align*}
		\end{example}

		\begin{example}
			Let $V$ be a Pareto random variable with cdf $F(x)=1-x^{-\alpha},\  \alpha >0, \ x>1 $. Let $w(x)=x^m,\ m>0, \ x>0 $, then we get
			\begin{align*}
				\Lambda_V^{w} (u) 
				&= w(F^{-1}(u)) f(F^{-1}(u))\\
				&=\alpha (1-u)^{\frac{\alpha - m+1}{\alpha}}.		
			\end{align*}
			Then from (\ref{minrssuformula}) we have 
			\begin{align*}
				J^{w}(\textbf{X}_{minRSSU}^{(n)}) &=-\frac{(n!)^2}{2}\prod_{i=1}^{n}\left(\int_{0}^{1} (1-u)^{2i-2} w(F^{-1}(u)) f(F^{-1}(u))du\right)\nonumber \\
				&=-\frac{(n!)^2}{2}\prod_{i=1}^{n}\left(\int_{0}^{1} (1-u)^{2i-2} \alpha (1-u)^{\frac{\alpha - m+1}{\alpha}}du\right)\nonumber \\
				&=-\frac{(n!)^2 \alpha^n}{2}\prod_{i=1}^{n} \frac{1}{2i\alpha -m+1}.
			\end{align*}
			
			Also, from (\ref{maxrssuformula}) we have
			\begin{align*}
				J^{w}(\textbf{X}_{maxRSSU}^{(n)}) &=-\frac{(n!)^2}{2}\prod_{i=1}^{n}\left(\int_{0}^{1} u^{2i-2} w(F^{-1}(u)) f(F^{-1}(u))du\right)\nonumber \\
				&=-\frac{(n!)^2}{2}\prod_{i=1}^{n}\left(\int_{0}^{1} u^{2i-2} \alpha (1-u)^{\frac{\alpha - m+1}{\alpha}}du\right)\nonumber \\	
				&=	-\frac{(n!)^2 \alpha^n}{2}\prod_{i=1}^{n} \beta\left(2i-1, \frac{2\alpha-m+1}{\alpha}\right).
			\end{align*}
			
		\end{example}

		\begin{example}\label{eg4uniform}
			Let $V$ be a uniform random variable with cdf $F(x)=x,\  \alpha >0, \ 0<x<1 $. Let $w(x)=x^m,\ m>0, \ x>0 $, then we get
			\begin{align*}
				\Lambda_V^{w} (u) 
				&= w(F^{-1}(u)) f(F^{-1}(u))=u^m.		
			\end{align*}
			Then from (\ref{minrssuformula}) we have 
			\begin{align*}
				J^{w}(\textbf{X}_{minRSSU}^{(n)}) &=-\frac{(n!)^2}{2}\prod_{i=1}^{n}\left(\int_{0}^{1} (1-u)^{2i-2} w(F^{-1}(u)) f(F^{-1}(u))du\right)\nonumber \\
				&=-\frac{(n!)^2}{2}\prod_{i=1}^{n}\left(\int_{0}^{1} (1-u)^{2i-2} u^mdu\right)\nonumber \\
				&=-\frac{(n!)^2 }{2}\prod_{i=1}^{n} \beta\left(2m+1,2i-1\right).
			\end{align*}
			
			Also, from (\ref{maxrssuformula}) we have
			\begin{align*}
				J^{w}(\textbf{X}_{maxRSSU}^{(n)}) &=-\frac{(n!)^2}{2}\prod_{i=1}^{n}\left(\int_{0}^{1} u^{2i-2} w(F^{-1}(u)) f(F^{-1}(u))du\right)\nonumber \\
				&=-\frac{(n!)^2}{2}\prod_{i=1}^{n}\left(\int_{0}^{1} u^{2i-2} u^m du\right)\nonumber \\	
				&=	-\frac{(n!)^2 }{2}\prod_{i=1}^{n} \frac{1}{m+2i-1}.
			\end{align*}	
		\end{example}
		
		\begin{example}
			Suppose that $X$ is a non-negative random variable from proportional hazard family (PHF) with cdf $F(x)=1-[\bar{F_0}(x)]^\theta $ and pdf $f(x)=\theta f_0(x)[\bar{F_0}(x)]^{\theta-1}$ where $F_0$ and $f_0$ stand for cdf and pdf of the underlying population, respectively, and $\theta>0$ is a shape parameter. Then  from (\ref{minrssuformulaline1}),
			\[J^{w}(\textbf{X}_{minRSSU}^{(n)}) =-\frac{\theta^{2n}(n!)^2}{2(2n\theta-1)!!}\prod_{i=1}^{n} E[w(F_0^{-1}(B_{1,2\theta i-1}))f_0(F_0^{-1}(B_{1,2\theta i-1}))].\]
		\end{example}
		
		\begin{example}
			Suppose that $X$ is a non-negative random variable from a proportional reversed hazard family (PRHF) with cdf $F(x)=[F_0(x)]^\theta $ and pdf $f(x)=\theta f_0(x)[F_0(x)]^{\theta-1}.$ Then from (\ref{maxrssuformulaline1}),
			\[J^{w}(\textbf{X}_{maxRSSU}^{(n)}) =-\frac{\theta^{2n}(n!)^2}{2(2n\theta-1)!!}\prod_{i=1}^{n} E[w(F_0^{-1}(B_{2\theta i-1,1}))f_0(F_0^{-1}(B_{2\theta i-1, 1}))]. \]
		\end{example}
		
		The conditions under which the GWJ will increase (decrease) are shown by the following Theorem. 
		\begin{theorem}\label{minrssuthm1}
			Let $X$ be a non-negative absolutely continuous non-negative random variable with pdf f and cdf F. Assume $\phi(x)$ is an increasing function and  $\frac{w(\phi(x))}{\phi^\prime (x)} \leq (\geq) w(x)$ and $\phi(0)=0$. If $Z=\phi(X)$, then 
			\begin{enumerate}[(i)]
				\item $J^{w}(\textbf{X}_{minRSSU}^{(n)})\leq (\geq) J^{w}(\textbf{Z}_{minRSSU}^{(n)})$
				\item $J^{w}(\textbf{X}_{maxRSSU}^{(n)})\leq (\geq) J^{w}(\textbf{Z}_{maxRSSU}^{(n)})$.
			\end{enumerate}
			
		\end{theorem}

		\noindent \textbf{Proof} (i)
		Let pdf and cdf of $Z$ be $h$ and $H$, respectively. Then 
		\begin{align*}
			\Lambda_Z^{w}(u)&=w\left(H^{-1}(u)\right) h\left(H^{-1}(u)\right)\\
			&=\frac{w(\phi(F^{-1}(u)))}{\phi^\prime(F^{-1}(u))}f(F^{-1}(u)) \hspace{5mm} \forall \hspace{5mm} 0<u<1.
		\end{align*}
		Note that $\phi(x)\geq \phi(0)$, $\forall$ $x\geq 0$. Hence for $ 0<u<1$, we have
		\begin{align*}
			\Lambda_Z^{w}(u)=\frac{w(\phi(F^{-1}(u)))}{\phi^\prime(F^{-1}(u))}f(F^{-1}(u))\leq w(F^{-1}(u))f(F^{-1}(u))=\Lambda_X^{w}(u).
		\end{align*}
		Therefore $J^{w}(\textbf{X}_{minRSSU}^{(n)})\leq J^{w}(\textbf{Z}_{minRSSU}^{(n)})$ using equation (\ref{minrssuformula}). 
		Proof of other parts can be done in a similar fashion.
		
		(ii) Similar to part (i). \hfill $\blacksquare$

		\section{Stochastic Comparision}
		
		In the following result, we provide the conditions for comparing two RSS schemes under different weights.
		
		\begin{theorem}\label{thm com rss1}
			Let $X$  and $Y$ be non negative random variables with pdf's $f$ and $g$, cdf's $F$ and $G$, respectively having $u_X=u_Y<\infty$.\\
			(a) If $w_1$ is increasing, $w_1(x)\geq w_2(x)$ and $X\le_{disp} Y$, then $J^{w_1}(\textbf{X}_{minRSSU}^{(n)})\le J^{w_2}(\textbf{Y}_{minRSSU}^{(n)})$.\\
			(b)   If $w_1$ is increasing, $w_1(x)\leq w_2(x)$ and $X\ge_{disp} Y$, then $J^{w_1}(\textbf{X}_{minRSSU}^{(n)})\ge J^{w_2}(\textbf{Y}_{minRSSU}^{(n)})$.\\
			(c) If $w_1$ is increasing, $w_1(x)\geq w_2(x)$ and $X\le_{disp} Y$, then $J^{w_1}(\textbf{X}_{minRSSU}^{(n)})\le J^{w_2}(\textbf{Y}_{minRSSU}^{(n)})$.\\
			(d)   If $w_1$ is increasing, $w_1(x)\leq w_2(x)$ and $X\ge_{disp} Y$, then $J^{w_1}(\textbf{X}_{maxRSSU}^{(n)})\ge J^{w_2}(\textbf{Y}_{maxRSSU}^{(n)})$.\\
		\end{theorem}
		\noindent \textbf{Proof}
		(a) Since $X\le_{disp} Y$, therefore we have $f(F^{-1}(u))\ge g(G^{-1}(u))$ for all $u\in (0,1)$. Then using Theorem 3.B.13(b) of Shaked and Shanthikumar (2007), $X\le_{disp} Y$ implies that $X\ge_{st} Y$. Hence $F^{-1}(u) \ge G^{-1}(u)$ $\forall$ $u\in (0,1)$. Since $w_1$ is increasing and $w_1(x)\geq w_2(x)$, then $w_1(F^{-1}(u)) \ge w_1(G^{-1}(u))\ge w_2(G^{-1}(u))$. 
		Hence 
		\begin{align}\label{stor1}
			\Lambda_X^{w_1} (u)&=w_1(F^{-1}(u)) f(F^{-1}(u))\nonumber\\
			&\ge w_2(G^{-1}(u)) g(G^{-1}(u))\nonumber\\
			=&\Lambda_Y^{w_2} (u).
		\end{align}
		Now, using (\ref{minrssuformula}),
		\begin{align*}
			J^{w_1}(\textbf{X}_{minRSSU}^{(n)})
			&=-\frac{(n!)^2}{2(n-1)!!}\prod_{i=1}^{n} E(\Lambda_X^{w_1}(B_{1,2i-1}))\\
			&\le -\frac{(n!)^2}{2(n-1)!!}\prod_{i=1}^{n} E(\Lambda_Y^{w_2}(B_{1,2i-1}))\\
			&=J^{w_2}(\textbf{Y}_{minRSSU}^{(n)}).\\
		\end{align*}
		(b) Proof is similar to part (a) \\
		(c) Proof is similar to part (a)\\
		(d) Proof is similar to part (a).\hfill $\blacksquare$

		If we take $w_1(x)=w_2(x)=w(x)$ in the above theorem, then The following corollary follows.
		
		\begin{corollary}\label{cor com rss1}
			Let $X$  and $Y$ be nonnegative random variables with pdf's $f$ and $g$, cdf's $F$ and $G$, respectively having $u_X=u_Y<\infty$; let $w$ be increasing. Then\\
			(a) If $X\le_{disp} Y$, then $J^{w}(\textbf{X}_{minRSSU}^{(n)})\le J^{w}(\textbf{Y}_{minRSSU}^{(n)})$.\\
			(b)  If $X\ge_{disp} Y$, then $J^{w}(\textbf{X}_{minRSSU}^{(n)})\ge J^{w}(\textbf{Y}_{minRSSU}^{(n)})$.\\
			(c) If $X\le_{disp} Y$, then $J^{w}(\textbf{X}_{maxRSSU}^{(n)})\le J^{w}(\textbf{Y}_{maxRSSU}^{(n)})$.\\
			(d)  If $X\ge_{disp} Y$, then $J^{w}(\textbf{X}_{maxRSSU}^{(n)})\ge J^{w}(\textbf{Y}_{maxRSSU}^{(n)})$.
		\end{corollary}

		\begin{remark}
			For $w_1(x)=1$ , the result in Corollary  \ref{cor com rss1} is proved by Qiu and Eftekharian (2021).
		\end{remark}
		
		Random variable $X$ is said to be smaller than $Y$ in the dispersive order, denoted $X\le_{disp} Y$, if $G^{-1}(F(x))-x$ is increasing in $x$ (Shaked and Shanthikumar (2007)).
		
		\begin{theorem}
			Let $X$  and $Y$ be nonnegative random variables with pdf's $f$ and $g$, cdf's $F$ and $G$, respectively. Let $\phi$ be a non-negative function such that $\phi^\prime(x)\geq 1$ for all $x.$ If $Y=\phi(X)$ and $w$ be increasing then $J^{w}(\textbf{X}_{minRSSU}^{(n)})\le J^{w}(\textbf{Y}_{minRSSU}^{(n)})$ and $J^{w}(\textbf{X}_{maxRSSU}^{(n)})\le J^{w}(\textbf{Y}_{maxRSSU}^{(n)})$.
		\end{theorem} 
		\noindent \textbf{Proof} Since $G^{-1}(F(x))-x= \phi(x)-1,$ therefore $X\le_{disp} Y$. Hence, the proof is completed from Corollary  \ref{cor com rss1}. \hfill $\blacksquare$
		
		\begin{corollary}\label{cor com rss2}
			Let $X$  and $Y$ be nonnegative random variables with pdf's $f$ and $g$, cdf's $F$ and $G$, respectively. If $Y=aX+b$ for $a\geq1, \ b \in \mathbb{R}$ and $w$ be increasing then $J^{w}(\textbf{X}_{minRSSU}^{(n)})\le J^{w}(\textbf{Y}_{minRSSU}^{(n)})$ and $J^{w}(\textbf{X}_{maxRSSU}^{(n)})\le J^{w}(\textbf{Y}_{maxRSSU}^{(n)})$.
		\end{corollary} 
		\noindent \textbf{Proof} Consider $\phi(x)=ax+b$ for $a\geq1, \ b \in \mathbb{R}.$ Hence, proof is completed from Corollary  \ref{cor com rss2}. \hfill $\blacksquare$

		\begin{lemma}\label{lemma1} [Ahmed et al. (1986)]
			Let $X$ and $Y$  be nonnegative random variables with pdf's $f$ and $g$, respectively, satisfying $f(0)\ge g(0)>0$. If $X\le_{su}Y$ (or $X\le_{*}Y$ or $X\le_{c}Y$), then $X\le_{disp}Y$.
		\end{lemma}
		
		One may refer to Shaked and Shanthikumar (2007) for details of convex transform order ($\leq_c$), star order ($\leq_{\star}$), super additive order ($\leq_{su}$), and dispersive order ($\leq_{disp}$). In view of Theorem \ref{thm com rss1}  and Lemma \ref{lemma1}, the following result is obtained.
		\begin{corollary}\label{corrolary3}
			Let $X$  and $Y$ be nonnegative random variables with pdf's $f$ and $g$, cdf's $F$ and $G$, respectively having $u_X=u_Y<\infty$ such that $f(0)\ge g(0)>0$.\\
			(a)  If $w_1$ is increasing, $w_1(x)\geq w_2(x)$ and $X\le_{su} Y$ (or $X\le_{*}Y$ or $X\le_{c}Y$), then $J^{w_1}(\textbf{X}_{minRSSU}^{(n)})\le J^{w_2}(\textbf{Y}_{minRSSU}^{(n)})$.\\
			(b)   If $w_1$ is increasing, $w_1(x)\leq w_2(x)$ and $X\ge_{su} Y$  (or $X\ge_{*}Y$ or $X\ge_{c}Y$), then $J^{w_1}(\textbf{X}_{minRSSU}^{(n)})\ge J^{w_2}(\textbf{Y}_{minRSSU}^{(n)})$.\\
			(c)  If $w_1$ is increasing, $w_1(x)\geq w_2(x)$ and $X\le_{su} Y$ (or $X\le_{*}Y$ or $X\le_{c}Y$), then $J^{w_1}(\textbf{X}_{maxRSSU}^{(n)})\le J^{w_2}(\textbf{Y}_{maxRSSU}^{(n)})$.\\
			(d)   If $w_1$ is increasing, $w_1(x)\leq w_2(x)$ and $X\ge_{su} Y$  (or $X\ge_{*}Y$ or $X\ge_{c}Y$), then $J^{w_1}(\textbf{X}_{maxRSSU}^{(n)})\ge J^{w_2}(\textbf{Y}_{maxRSSU}^{(n)})$.
		\end{corollary}
		
		Let X be a nonnegative random variable with an absolutely
		continuous distribution function. Then X has a decreasing density if, and only if, $U \leq_c X,$ where U is a uniform(0, 1) random variable (Shaked and Shanthikumar (2007), Example 4.B.12). From Example \ref{eg4uniform}, 	$J^{w}(\textbf{X}_{minRSSU}^{(n)}) =-\frac{(n!)^2 }{2}\prod_{i=1}^{n} \beta\left(2m+1,2i-1\right)$	and 	   $J^{w}(\textbf{X}_{maxRSSU}^{(n)}) =	-\frac{(n!)^2 }{2}\prod_{i=1}^{n} \frac{1}{m+2i-1}.$ Thus, we get the lower bound of GWJ of minRSSU and maxRSSU data by Corollary \ref{corrolary3} and Example 4.B.12 of Shaked and Shanthikumar (2007) when considering $w_1(x)=w_2(x)=x^m, \ m\geq0.$ 	
		\begin{proposition}
			Let X be a non-negative random variable with decreasing pdf $f$ such that $f(0) \leq 1$ and $w(x)=x^m, \ m\geq 0 $ then $J^{w}(\textbf{X}_{minRSSU}^{(n)}) \geq -\frac{(n!)^2 }{2}\prod_{i=1}^{n} \beta\left(2m+1,2i-1\right)$	and 	   $J^{w}(\textbf{X}_{maxRSSU}^{(n)}) \geq	-\frac{(n!)^2 }{2}\prod_{i=1}^{n} \frac{1}{m+2i-1}.$
		\end{proposition}

		\begin{lemma}\label{lemmashaked4.B.11}
			(Theorem 4.B.11, Shaked and Shanthikumar (2007)) Let Y denote the exponential random variable with parameter $\lambda$ and mean $\frac{1}{\lambda}$. Let X be a nonnegative random variable. Then \\
			(i) X is IFR if and only if $X \leq_{c}  Y$, \\
			(ii) X is IFRA if and only if $X \leq_{*} Y$, and \\
			(iii) X is NBU if and only if $ X\leq_{su} Y$.
		\end{lemma}

		We obtained $J^{w}(\textbf{Y}_{minRSSU}^{(n)})=-\frac{1}{2}\left(\frac{1}{\lambda^{m-1}}\right)^{n}\left(\frac{\Gamma(m+1)}{2^{m+1}}\right)^n \left(\frac{1}{n!}\right)^{m-1}$ in Example \ref{eg2exp} where $Y$ denote exponential random variable with parameter $\lambda$ and mean $\frac{1}{\lambda}$. We can write the upper bound of GWJ of minRSSU data by Example \ref{eg2exp}, Corollary \ref{corrolary3} and Lemma \ref{lemmashaked4.B.11} when considering $w_1(x)=w_2(x)=x^m, \ m\geq0.$		
		\begin{proposition}
			If $X$ is IFR (IFRA, NBU), $f(0) \geq \lambda$ and $w(x)=x^m, \ m\geq 0 $, then \[J^{w}(\textbf{X}_{minRSSU}^{(n)})\leq-\frac{1}{2}\left(\frac{1}{\lambda^{m-1}}\right)^{n}\left(\frac{\Gamma(m+1)}{2^{m+1}}\right)^n \left(\frac{1}{n!}\right)^{m-1}.\]
		\end{proposition}
		
		\begin{lemma}\label{lemma3}
			(Bagai and Kochar (1986)) Let $X$ and $Y$ be two non-negative random variables. If $X\geq_{hr} (\leq_{hr}) Y$
			and either $X$ or $Y$ is DFR (decreasing failure rate), then $X\geq_{disp} (\leq_{disp}) Y .$
		\end{lemma}
		
		\begin{lemma}\label{3.b.28shaked}
			(Theorem 3.B.28 Shaked and Shanthikumar (2007))	Let $X_{i:n}$ and $X_{j:m}$ denote the ith and jth order statistics from a DFR distribution $F$ of size $n$ and $m,$ respectively. Then $X_{i:n} \leq_{disp} X_{j:m}$  for $i \leq j$ and $n-i\geq m-j.$
		\end{lemma}
		One more application of Corollary \ref{cor com rss1} is illustrated in the following proposition which provides us lower bound of $J^{w}(\textbf{X}_{maxRSSU}^{(n)}).$
		
		\begin{proposition}
			Let $X$ be a non negative random variable with $r_X(x)\leq \lambda$ for all $x \geq 0$ and $w(x)=x^m, \ m\geq 0$ then \[J^{w}(\textbf{X}_{maxRSSU}^{(n)})\geq  -\frac{1}{2} \left(\frac{\lambda^2 \Gamma(m+1)}{(2\lambda)^{m+1}}\right)^n. \]
		\end{proposition}
		\noindent \textbf{Proof} Let $Z$ have the exponential distribution with failure rate $\lambda$.  Since, $r_X(x)\leq \lambda =r_Z(x)$ and $Z$ has DFR property,  therefore using Lemma \ref{lemma3}, $X \geq_{disp} Z.$ Also using Lemma \ref{3.b.28shaked}, $Z_{1:1}\leq_{disp} Z_{2:2} \leq_{disp} Z_{3:3} \leq_{disp} \dots \leq_{disp} Z_{n:n},$ which leads to $ J^w(Z_{1:1})\leq J^w(Z_{2:2}) \leq J^w(Z_{3:3}) \dots J^w(Z_{n:n}),$ because of  Corollary \ref{cor com rss1}. Using Corollary \ref{cor com rss1}, we get 
		
		\begin{align*}
			J^{w}(\textbf{X}_{minRSSU}^{(n)}) &\geq J^{w}(\textbf{Z}_{minRSSU}^{(n)}) = -\frac{1}{2} \prod_{i=1}^{n} [-2J^w(Z_{i:i})] \\
			&\geq -\frac{1}{2} [-2J^w(Z_{1:1})]^n = -\frac{1}{2} \left(\frac{\lambda^2 \Gamma(m+1)}{(2\lambda)^{m+1}}\right)^n. 
		\end{align*}
		Hence, we obtain the required result. \hfill $\blacksquare$

		One may ask that whether the condition $X\le_{disp}Y$ in Theorem \ref{thm com rss1}  may be relaxed by $J^{w_1}(X)\leq J^{w_2}(Y)$. The following result gives a positive answer to this assertion. Proof of the following theorem is similar to the proof of Theorem 5.3 of Gupta and Chaudhary (2023), hence omitted.		
		\begin{theorem}\label{thm com rss3}
			Let $X$  and $Y$ be nonnegative random variables with pdf's $f$ and $g$, cdf's $F$ and $G$, respectively. Let $\Delta (u)=w_1(F^{-1}(u))f(F^{-1}(u))-w_2(G^{-1}(u))g(G^{-1}(u))$,
			\[A_1=\{0\le u \le 1|\Delta (u)>0\},\ A_2=\{0\le u \le 1|\Delta (u)<0\}.\]
			If $ \inf_{A_1}\phi_{2i-1:2n-2i}(u)\geq  \sup_{A_2}\phi_{2i-1:2n-2i}(u)$, and if $J^{w_1}(X)\leq J^{w_2}(Y)$, then  $J^{w_1}(\textbf{X}_{minRSSU}^{(n)})$ $\le J^{w_2}(\textbf{Y}_{minRSSU}^{(n)}).$
		\end{theorem}
		
		If we take $w_1(x)=w_2(x)$ in the above theorem, then we have the following corollary.
		\begin{corollary}\label{cor com rss3}
			Let $X$  and $Y$ be nonnegative random variables with pdf's $f$ and $g$, cdf's $F$ and $G$. Let $\Delta (u)=w_1(F^{-1}(u))f(F^{-1}(u))-w_1(G^{-1}(u))g(G^{-1}(u))$,
			\[A_1=\{0\le u \le 1|\Delta (u)>0\},\ A_2=\{0\le u \le 1|\Delta (u)<0\}.\]
			If  $\ \inf_{A_1}\phi_{2i-1:2n-2i}(u)\geq  \sup_{A_2}\phi_{2i-1:2n-2i}(u)$, and if $J^{w_1}(X)\leq J^{w_1}(Y)$, then  $J^{w_1}(\textbf{X}_{minRSSU}^{(n)})$ $\le J^{w_1}(\textbf{Y}_{minRSSU}^{(n)}).$
		\end{corollary}

		\section{Monotone Properties}
		In this section, we give a condition under which the general weighted extropy of minRSSU data is decreasing in the sample size n. 
		
		\begin{theorem}
			Let $\textbf{X}_{minRSSU}^{(n)}$ be a minRSSU sample from a distribution with cdf $F$ and pdf $f.$ If $w(F^{-1}(u)) f(F^{-1}(u)) \geq 1$ for all $0<u<1,$ then $J^w(\textbf{X}_{minRSSU}^{(n)})$ is decreasing in $n\geq 1.$
		\end{theorem}
		\noindent \textbf{Proof} Using (\ref{minrssuformulaline1}), it follows that
		\begin{align*}
			\frac{J^w(\textbf{X}_{minRSSU}^{(n+1)})}{J^w(\textbf{X}_{minRSSU}^{(n)})}&=-2J^w(X_{1:n+1})\\
			&=\int_{0}^{+\infty} w(x) f^2_{1:n+1}(x)dx \\
			&= (n+1)^2 \int_{0}^{1} w(F^{-1}(u)) (1-u)^{2n} f(F^{-1}(u))du \\
			& \geq \frac{(n+1)^2}{2n+1} \\
			& \geq 1.
		\end{align*}
		Therefore, $J^w(\textbf{X}_{minRSSU}^{(n+1)}) \leq J^w(\textbf{X}_{minRSSU}^{(n)})$ since both are non-positive. Hence, $J^w(\textbf{X}_{minRSSU}^{(n)})$ is decreasing in $n\geq 1.$ \hfill $\blacksquare$
		
		The general weighted extropy of the maxRSSU data is also decreasing in the sample size $n$ like the above theorem.
		\begin{theorem}
			Let $\textbf{X}_{maxRSSU}^{(n)}$ be a maxRSSU sample from a distribution with cdf $F$ and pdf $f.$ If $w(F^{-1}(u)) f(F^{-1}(u)) \geq 1$ for all $0<u<1,$ then $J^w(\textbf{X}_{maxRSSU}^{(n)})$ is decreasing in $n\geq 1.$
		\end{theorem}
		\noindent \textbf{Proof} Using (\ref{maxrssuformulaline1}), it follows that
		\begin{align*}
			\frac{J^w(\textbf{X}_{maxRSSU}^{(n+1)})}{J^w(\textbf{X}_{maxRSSU}^{(n)})}&=-2J^w(X_{n+1:n+1})\\
			&=\int_{0}^{+\infty} w(x) f^2_{1:n+1}(x)dx \\
			&= (n+1)^2 \int_{0}^{1} w(F^{-1}(u)) u^{2n} f(F^{-1}(u))du \\
			& \geq \frac{(n+1)^2}{2n+1} \\
			& \geq 1.
		\end{align*}
		Therefore, $J^w(\textbf{X}_{maxRSSU}^{(n+1)}) \leq J^w(\textbf{X}_{maxRSSU}^{(n)})$ since both are non-positive. Hence, $J^w(\textbf{X}_{maxRSSU}^{(n)})$ is decreasing in $n\geq 1.$ \hfill $\blacksquare$

		\begin{corollary}
			If $X$ is PHF with $\theta >1/4 \ \text{and}\ w(F_0^{-1}(u)) f_0(F_0^{-1}(u)) \geq 1$ for all $0<u<1.$ Then, $J^w(\textbf{X}_{minRSSU}^{(n)})$ is decreasing in $n\geq 1.$ 
		\end{corollary}

		\begin{corollary}
			If $X$ is PRHF with $\theta >1/4 \ \text{and}\ w(F_0^{-1}(u)) f_0(F_0^{-1}(u)) \geq 1$ for all $0<u<1.$ Then, $J^w(\textbf{X}_{maxRSSU}^{(n)})$ is decreasing in $n\geq 1.$ 
		\end{corollary}

		\section{Bounds}
		In this section, we obtain bounds for $J^w(\textbf{X}_{minRSSU}^{(n)})$ and $J^w(\textbf{X}_{maxRSSU}^{(n)})$ in terms of $J^w(\textbf{X}_{SRS}^{(n)}).$
		
		\begin{theorem} If $\textbf{X}_{minRSSU}^{(n)}, \textbf{X}_{maxRSSU}^{(n)} \text{and} \ \textbf{X}_{SRS}^{(n)}$ denote minRSSU sample, maxRSSU sample and SRS sample respectively. Then we have 
			\begin{enumerate}[(i)]
				\item  $J^w(\textbf{X}_{minRSSU}^{(n)}) \geq (n!)^2 J^w(\textbf{X}_{SRS}^{(n)}).$ 
				\item  $J^w(\textbf{X}_{maxRSSU}^{(n)}) \geq (n!)^2 J^w(\textbf{X}_{SRS}^{(n)}).$
			\end{enumerate}	
		\end{theorem}
		\noindent \textbf{Proof} (i) From (\ref{minrssuformula}),
		\begin{align*}
			J^w(\textbf{X}_{minRSSU}^{(n)}) 
			&= -\frac{(n!)^2}{2}\prod_{i=1}^{n}\left(\int_{0}^{1}w(F^{-1}(u)) (1-u)^{2i-2}f(F^{-1}(u))du\right) \\
			& \geq -\frac{(n!)^2}{2}\prod_{i=1}^{n}\left(\int_{0}^{1}w(F^{-1}(u)) f(F^{-1}(u))du\right)\\
			&= (n!)^2 J^w(\textbf{X}_{SRS}^{(n)}).
		\end{align*}
		(ii) Proof follows in a similar fashion to part (i). \hfill $\blacksquare$
		
		We denote $X_{i:n}$ as the $i$th order statistic based on a simple random sample from $X$ of size $n,$ then under the assumption of perfect judgment ranking, the ranked set sample is defined as $X^{(n)}_{RSS}=\{X_{i:n}, i=1,\dots,n\}.$ Note that $X^{(n)}_{RSS}$ is a collection of order statistics from $X$ but independent, differing from the order statistics of a simple random sample. Thus, the general weighted extropy of $X^{(n)}_{RSS}$
		was defined by Gupta and Chaudhary (2023) as

		\begin{align}\label{GWJofRSS}
			J^{w}(\textbf{X}_{RSS}^{(n)})=-\frac{1}{2}\prod_{i=1}^{n}\left(-2J^{w}(X_{(i:n)})\right).
		\end{align}

		Next, we will demonstrate that if the underlying distribution has the DFR ageing property, the general weighted extropy of maxRSSU data is higher than that of RSS and MinRSSU data.
		\begin{theorem}\label{boundDFRminmaxrss}
			If $\textbf{X}_{minRSSU}^{(n)}, \textbf{X}_{maxRSSU}^{(n)} \text{and} \ \textbf{X}_{SRS}^{(n)}$ denote minRSSU , maxRSSU  and SRS  data from distribution $X$ with DFR ageing property, respectively. Let $w(x)$ be an increasing function. Then we have 		
			\[J^w(\textbf{X}_{maxRSSU}^{(n)}) \geq  J^w(\textbf{X}_{RSS}^{(n)}) \geq J^w(\textbf{X}_{minRSSU}^{(n)}). \] 
		\end{theorem}	
		\noindent \textbf{Proof} It follows from Lemma \ref{3.b.28shaked} that if $X_1,\dots, X_n, \dots$ is a sequence of independent and identically distributed observations from DFR distribution, then $X_{1:i}\leq_{disp} X_{i:n}\leq_{disp}X_{i:i},$ for all $i=1,2,\dots, n,$ which implies $J^w(X_{1:i})\leq J^w(X_{i:n})\leq J^w(X_{i:i})$ by Corollary 2.1 in Gupta and Chaudhary (2023). From (\ref{minrssuformulaline1}), (\ref{maxrssuformulaline1}) and (\ref{GWJofRSS}), result follows. \hfill $\blacksquare$
		
		\begin{example}
			Let $Z$ have an exponential distribution with cdf $F_Z(z)=1-e^{-\lambda z}, \ \lambda >0, \ z>0$. Since $Z$ has DFR property, therefore using Theorem \ref{boundDFRminmaxrss}, we obtain $J^w(\textbf{Z}_{maxRSSU}^{(n)}) \geq  J^w(\textbf{Z}_{RSS}^{(n)}) \geq J^w(\textbf{Z}_{minRSSU}^{(n)}).$
		\end{example}
		In next theore, we find lower bound of $J^w(\textbf{Z}_{maxRSSU}^{(n)})$ and $J^w(\textbf{Z}_{minRSSU}^{(n)})$ when upper bounds of pdf and $w(x)$ are known.
		\begin{theorem}
			Let $X$ be a random variable with mode $m,$ and $f(m)=M$ i.e., $0\leq f(x)\leq M$  for all x and $0 \leq w(x)\leq k, \ k>0.$ Then for all $n\geq1$, we have,	
			$J^w(\textbf{X}_{minRSSU}^{(n)}) \geq - \frac{(n!)^2 k^n M^n}{2(2n-1)!!}$ and $J^w(\textbf{X}_{maxRSSU}^{(n)}) \geq - \frac{(n!)^2 k^n M^n}{2(2n-1)!!}.$
		\end{theorem}
		\noindent \textbf{Proof} From (\ref{minrssuformula}), we have
		\begin{align*}
			J^w(\textbf{X}_{minRSSU}^{(n)}) 
			&= -\frac{(n!)^2}{2}\prod_{i=1}^{n}\left(\int_{0}^{1}w(F^{-1}(u)) (1-u)^{2i-2}f(F^{-1}(u))du\right) \\
			& \geq -\frac{(n!)^2}{2}\prod_{i=1}^{n}\left(\int_{0}^{1} kM (1-u)^{2i-2}du\right)\\
			&= -\frac{(n!)^2 k^n M^n}{2(2n-1)!!}.
		\end{align*}
		
		Similarly using (\ref{maxrssuformula}), proof of the other part follows. \hfill $\blacksquare$
		
		Now,  we derive a lower bound for the general weighted extropy of  MinRSSU and maxRSSU data in terms of known expectations of $w(X)f^2(X).$
		
		\begin{theorem}
			Let $X$ be a non-negative random variable with pdf $f.$ Let $\kappa_f=E[(w(X)f(X))^2].$ Then, 
			$J^w(\textbf{X}_{maxRSSU}^{(n)}) \geq - \frac{(n!)^2 (\kappa_f^{n/2})}{2 \sqrt{\prod_{i=1}^{n} (4i-3)}}$ and $J^w(\textbf{X}_{minRSSU}^{(n)}) \geq - \frac{(n!)^2 (\kappa_f^{n/2})}{2 \sqrt{\prod_{i=1}^{n} (4i-3)}}.$ 	
		\end{theorem}
		\noindent \textbf{Proof} From (\ref{maxrssuformula}) and by using Cauchy Schwarz inequality, we have, 
		\begin{align*}
			J^w(\textbf{X}_{maxRSSU}^{(n)}) &=-\frac{(n!)^2}{2(2n-1)!!}\prod_{i=1}^{n}\left(\int_{0}^{1} (2i-1)u^{2i-2}w(F^{-1}(u)) f(F^{-1}(u))du\right) \\
			&\geq -\frac{(n!)^2}{2(2n-1)!!}\prod_{i=1}^{n} \left[\int_{0}^{1} (2i-1)^2u^{4i-4}du \right]^{\frac{1}{2}} \left[\int_{0}^{1}w^2(F^{-1}(u)) f^2(F^{-1}(u))du\right]^{\frac{1}{2}} \\
			&= -\frac{(n!)^2 \kappa_f^{n/2}}{2(2n-1)!!}\prod_{i=1}^{n} \left[\int_{0}^{1} (2i-1)^2u^{4i-4}du \right]^{\frac{1}{2}} \\
			&= -\frac{(n!)^2 \kappa_f^{n/2}}{2(2n-1)!!}\prod_{i=1}^{n} \left[\frac{(2i-1)^2}{4i-3}\right]^{\frac{1}{2}}  	\\
			&= -\frac{(n!)^2 \kappa_f^{n/2}}{2 \sqrt{\prod_{i=1}^{n}(4i-3)}}.  	
		\end{align*}
		Similarly from (\ref{minrssuformula}), proof of other parts follows. \hfill $\blacksquare$  
		
		\section{Conclusion} \label{s10conclusion}\label{section8}
		In this work, we considered the GWJ of ranked set sampling with unequal sample size. We obtained expressions for GWJ of minRSSU and maxRSSU. Under some assumption on weight $w(x),$ GWJ value comparison of minRSSU and maxRSSU have been done. The behaviour of GWJ of minRSSU and maxRSSU have been given under some stochastic order assumptions. Bounds on GWJ of GWJ of minRSSU and maxRSSU have been derived due to the DFR property of exponential distribution. We have shown that the GWJ of minRSSU and maxRSSU are decreasing with respect to sample size under certain assumptions. Bounds on GWJ of minRSSU, maxRSSU, RSS and SRS have been obtained in terms of one another.\\
		\\

		\noindent \textbf{\Large Funding} \\
		\\
		Santosh Kumar Chaudhary would like to thank the Council of Scientific and industrial research (CSIR), Government of India (File Number 09/0081 (14002)/2022-EMR-I) for financial assistance.\\
		\\
		\textbf{ \Large Conflict of interest} \\
		\\
		The authors declare no conflict of interest.

	\end{document}